\documentclass[12pt,a4paper]{amsart}
\newtheorem{theorem}{Theorem}[section]

\newtheorem{lemma}[theorem]{Lemma}

\newtheorem{remark}[theorem]{Remark}

\newtheorem{notation}[theorem]{Notation}
\newtheorem{definition}[theorem]{Definition}

\newcommand{\N}{\mathbb{N}}

\title{On the frequency of permutations containing a long cycle} 
\author{Alice C. Niemeyer and Cheryl E. Praeger}
\date{\today}

\begin{document}
\begin{abstract}
A general explicit upper bound is obtained for the proportion $P(n,m)$
of elements of  order dividing $m$, where $n-1 \le m  \le cn$ for some
constant $c$,  in the finite symmetric  group $S_n$.  This  is used to
find lower bounds for the conditional probabilities that an element of
$S_n$  or $A_n$  contains an  $r$-cycle,  given that  it satisfies  an
equation  of the  form $x^{rs}=1$  where $s\leq3$.   For  example, the
conditional  probability that an  element $x$  is an  $n$-cycle, given
that $x^n=1$, is always greater  than $2/7$, and is greater than $1/2$
if $n$ does  not divide $24$.  Our results  improve estimates of these
conditional probabilities  in earlier work of the  authors with Beals,
Leedham-Green  and   Seress,  and  have   applications  for  analysing
black-box  recognition   algorithms  for  the   finite  symmetric  and
alternating groups.
\end{abstract}

\maketitle
 \section{Introduction}
This paper  is concerned with  the proportions of permutations  in the
symmetric and alternating groups  on $n$ points satisfying an equation
$x^m = 1$ for  various values of $m = O(n)$, that is,  for $m \le k n$
for some  constant $k$. Our  interest in such permutations  stems from
their use in `black-box  algorithms' to recognise finite symmetric and
alternating  groups. 
We give  a detailed
analysis of the proportion of  elements in $S_n$ of order dividing $m$
for large values of $n$, where  $m=O(n)$ and $m\geq n$
in \cite{NiemeyerPraeger05b}. These results 
 provide    asymptotic   bounds   for   $n$
sufficiently  large.  The  focus  of this  paper  is to find  explicit
probability  bounds  for  all   $n$.  
A shortened version of this paper
will be published \cite{NiemeyerPraeger05c} in the Journal of Algebra.
Such  bounds  are  useful in
algorithmic  applications: the  bounds obtained  in this  paper are
significant    improvements   on    explicit   estimates    given   in
\cite{Bealsetal03} for  the proportions of elements in  $S_n$ or $A_n$
satisfying   certain  equations,   and   the  associated   conditional
probabilities needed for the  algorithms.  
To explain the relevance of
these equations in  the design of the algorithms  we make some general
remarks about these black-box algorithms in Section~\ref{sec:appl}.

\subsection{Statement of results}
The purpose of this paper is first to prove in Theorem~\ref{first} a
general upper bound for the proportion $P(n,m)$ of elements in $S_n$ of order
dividing $m$, where $m=O(n)$. 
It has applications beyond those of this
paper, see~\cite{NiemeyerPraeger05}.
We then obtain in Theorem~\ref{melbounds_sn}
explicit lower bounds for the conditional
probabilities that an element of $S_n$ or $A_n$ has a relevant cycle
structure, given that it satisfies a certain equation. 

The  statement of Theorem~\ref{first} uses the following integer function:
\begin{definition}\label{def:gamma}
\begin{equation*}
\gamma(m) = 
  \begin{cases}
2  & \mbox{if\ } 360 < m,\\
 2.5& \mbox{if\ } 60 < m  \le 360,\\
 3.345 &\mbox{if\ } m \le 60.
	  \end{cases}
\end{equation*}
\end{definition}

\begin{theorem}\label{first}
Let $n,m$ be positive integers such that $m\geq n-1.$ Then the
proportion $P(n,m)$ of elements of $S_n$ of order dividing $m$ satisfies
\[
P(n,m)\leq  \frac{1}{n}+\frac{\gamma(m) m}{n^2}.
\]
\end{theorem}

If $m$ is very much larger than $n$ then the upper bound is greater
than 1 and hence of no use. However, if, say, $360 < m \le cn$ for
some constant $c$, then Theorem~\ref{first} implies that
$P(n,m) \le (2c+1)/n.$ 
It is difficult to give lower bounds for $P(m,n)$ that hold for all
$m$ and $n$.  
However, for example for a non-negative integer $k$, if $n-k$ divides
$m$ then $P(n,m) \ge \frac{1}{k!n}$, while if $n$ is even and
$n/2 -k$ divides $m$
then $P(n,m) \ge \frac{2}{(2k)!n^2}.$

Table~1   lists  the  kinds   of  elements   $g$  the   algorithms  in
\cite{Bealsetal03} seek  in $S_n$  or $A_n$ with  $n$ as given  in the
second column.  The fourth column headed {\sc CycType} lists the cycle
type of the  element $g$ in terms of a parameter  $r$ which is defined
in the third column. The fifth column records the order of $g^r$ 
and the last column records the group, either $A_n$ or $S_n$,
containing $g$.  Note  
that we omit  fixed points in the cycle notation.  Thus, for example a
permutation  in  $S_n$ with  cycle  type  $2^1(n-3)^1$  has one  fixed
point. In Section~\ref{sec:appl}  we give a brief account  of the role
of these elements in recognition algorithms for $S_n$ and $A_n$.

\bigskip
\begin{center}\label{tbl:cyctype}
\begin{tabular}{|c|c|lcl|c|}\hline
Case & $n$&$r$&{\sc CycType} & $|g^r|$ & $G$ \\ \hline
1 &    & $n$ & $r^1$  & 1 & $S_n$ \\
2 &  odd & $n-2$ & $2^1 r^1$ & 2 & $S_n$ \\
3 & even & $n-3$ & $2^1 r^1$ & 2 & $S_n$\\
4 &odd   & $n$ & $r^1$  & 1 &$A_n$ \\
5 & even   & $n-1$ & $r^1$ & 1 &$A_n$\\
6 & $2$ or $4 \pmod{6}$&$n-3$&$3^1r^1$& $3$ & $A_n$\\
7 & $3$ or $5 \pmod{6}$&$n-4$&$3^1r^1$ &$3$& $A_n$\\
8 & $0\pmod{6}$       &$n-5$&$3^1r^1$& $3$& $A_n$\\
9 & $1\pmod{6}$       &$n-6$&$3^1r^1$& $3$&$A_n$\\
10 & $1\pmod{6}$       &$n-5$&$2^13^1r^1$& $3$ &$A_n$\\ \hline
\end{tabular}

\medskip
{\bf Table~1\quad Relevant cycle types}
\end{center}

Theorem~\ref{melbounds_sn} gives
our improvements to the  estimates in \cite{Bealsetal03}
on conditional  probabilities for finding  such elements in $S_n$  or $A_n$.
Recall the integer function $\gamma(m)$ defined in Definition~\ref{def:gamma}.

\begin{theorem}\label{melbounds_sn}
Let $n\geq5$, let $G$ and $r$ be as in Table~$1$, and let $g$ be a uniformly
distributed random element of $G$. 
Then
\begin{enumerate}
\item[(a)] In cases $1, 4$ or $5$ of Table~$1$ let\\
${\displaystyle P =
{\rm Prob}\big( g  \ 
\mbox{is an $r$-cycle}\,|\,g^r=1\big)}.$ Then 
$$P \ge  1 - \frac{8+15\gamma(n)}{n^{2/3}}. \mbox{\ 
Also,\  } P \ge \begin{cases} 1/2 & \mbox{if\ } n \not\,\,\mid 24, \\
2/7 &  \mbox{if\ } n \mid 24. \\
		    \end{cases}
$$

\item[(b)] In cases $2$ or $3$ of Table~$1$ let $n\geq 8$ 
and\\ 
${\displaystyle P = 
{\rm Prob}\big( g  \ 
\mbox{has an $r$-cycle}\,|\,g^{2r}=1\ \mbox{and}\ |g^r|=2\big)}.$
Then 
$$P \ge 1 - \frac{18 + 76\gamma(2r)}{n^{2/3}}. \mbox{\ Also,\ }
P \ge \begin{cases} 1/3 & \mbox{if\ } n \not= 11, 17, 18,\\
1/4 &  \mbox{if\ } n = 11, 17, 18. \\
		    \end{cases}
$$
\item[(c)] In cases $6-10$ of Table~$1$ let $n\geq 8$  and \\
${\displaystyle P = {\rm Prob}\big( g \mbox{\ has\ an\ } r-\mbox{cycle\ }
\mid \,g^{3r}=1\ \mbox{and}\ |g^r|=3 
\big).}$  Then $$P \ge \begin{cases}
1 - \frac{98+839\gamma(3r)}{n^{2/3}}&\mbox{in\ cases\ $6-8, 10$ }\\
\frac{1}{2} - \frac{46+228\gamma(3r)}{n^{2/3}}&\mbox{in\ case\ $9$}.\\
			 \end{cases}
$$
Moreover, the lower bounds on $P$ given in Table~$2$ hold.
\begin{center}\label{tbl:two}
\begin{tabular}{|c|c|c|}\hline
P  $\ge$ & Case &$(n,r)$ \\ \hline
$3/10$ & $9$ & $(31, 25)$ \\
$3/10$ & $10$ & $(185, 80)$ \\
$3/20$ & $10$ & $(13, 8)$ or $(25, 20 )$\\
$1/3$ &  \multicolumn{2}{c|}{otherwise} \\
\hline
\end{tabular}\\
\medskip
{\bf Table~2\quad Lower bounds}
\end{center}
\end{enumerate}
\end{theorem}

Note  that  in \cite[Corollary~1.4]{NiemeyerPraeger05b}  we  prove
better  asymptotic 
bounds, namely  $1  - \frac{c}{n} +  O(\frac{1}{n^{1.5-o(1)}})$ 
for sufficiently  large $n$ (where the constant $c$ depends on the case),
for the  conditional  probabilities 
in Theorem~\ref{melbounds_sn} in the case where $G = S_n$.
However  these  bounds are  valid  only  for $n$  ``sufficiently
large'', whereas explicit lower bounds are required for each value of
$n$ for the algorithms. 
  
\subsection{Brief comments on our approach}

The key ingredient that enabled us to achieve our
results  was our two stage approach to the analysis. 
In the first stage we obtained in Theorem~\ref{first} a 
uniform upper bound  for the proportion $P(n,m)$ of permutations in
$S_n$
of order dividing $m$ for all $n, m$ with $n-1\leq m$. Although 
rather weak if $m$ is much larger than $n$ this result enabled us, in the 
second stage of our analysis, to obtain in Theorem~\ref{melbounds_sn}
explicit 
and improved bounds when applied to the families of 
permutations needed for the algorithms in \cite{Bealsetal03}.  One reason for 
this success was the precision we achieved in estimating the parameter 
$\gamma(m)$ in Theorem~\ref{first}. The idea behind the proof of 
Theorem~\ref{first} is a refinement of the approach of Beals and Seress 
in the proof of 
\cite[Theorem~3.7]{Bealsetal03} to study the cycles of elements of $S_n$ 
containing three specified points.

In Section~\ref{sect:divisors} we collect some well-known or sharpened
versions of  well-known upper bounds on  the number of  divisors of an
integer.   Section~\ref{sec:proportions}   is  devoted  to  elementary
properties  of  $P(n,m)$ and  the  proof  of Theorem~\ref{first}.   In
Section~\ref{sec:longcycle}  we derive the  practical upper  bounds on
the conditional probability that a  random element $g\in S_n$ or $A_n$
of order dividing a certain number 
has cycle type {\sc CycType} as specified in one
of the rows  of Table~1 (with $|g^r|$ as  in the last entry
of that row).

\subsection{Black-box algorithms for recognising $S_n$ and $A_n$}
\label{sec:appl} 
Black-box  algorithms  make  few  assumptions  about  how  groups  are
represented: group  elements may  be multiplied, inverted,  and tested
for equality. These three  operations are called black-box operations,
and no other  operations are permitted. Elements of  a black-box group
are represented as strings of zeros and ones, and the lengths of these
strings for a finite group $G$ can therefore be taken as approximately
$\log_2(|G|)$. Black-box algorithms are  then regarded as efficient if
the  number of  black-box  operations they  require  is polynomial  in
$\log_2(|G|)$,  that is,  at most  $O((\log|G|)^c)$ for  some constant
$c$. For  example, if $G$ is  the finite symmetric group  $S_n$ or the
alternating  group  $A_n$,  then  $\log|G|=O(n\log n)$,  so  efficient
black-box computations in these  groups should take $O(n^c)$ black-box
operations, for some constant $c$.

In particular in a black-box group $G = S_n$ it is
impossible  to determine
conclusively by black-box operations
 the cycle structure of an arbitrary element. Furthermore,
computing the order of a random element $g$ is too expensive since the
average  value for  $\log|g|$ is  $(1/2)\log^2n$, and  $g$  might have
order as  large as $e^{(1+o(1))(n\log n)^{1/2}}$  (see \cite{ETa, ETb}
and \cite[p. 222]{Landau09}). It  is however feasible to check whether
an element  satisfies an equation of  the form $x^m=1$  with less than
$2\log_2m$ black-box  operations (by the method  of repeated squaring)
and  this cost is acceptable  provided that  the number  of elements  to be
considered is not too great.

The $S_n$-recognition algorithm in \cite{Bealsetal03} has several
components and its analysis is based on the availability of
independent uniformly distributed random elements of the input group. 
It takes as input a black-box group $G$. If $G$ is isomorphic to $S_n$
then with high probability it returns an isomorphism $\lambda: G
\rightarrow S_n.$ 
In the first two steps elements $x$ and $y$ are constructed
such that $\lambda(x)$ is an $n$-cycle and $\lambda(y)$ is a 
transposition. Next, a random conjugate $y'$ of
$y$ is sought such that,
if $G$ is isomorphic to $S_n$, with high probability $x, y'$ satisfy a
standard presentation for $S_n$. Checking the presentation guarantees
that the subgroup $\langle x, y'\rangle$ generated by $x$ and $y'$ is
isomorphic to $S_n$. A further
algorithm tests that each generator of the input group $G$ lies
in the subgroup $\langle x, y'\rangle$ completing the proof that $G
\cong S_n$. A recognition algorithm for $A_n$ proceeds in a similar
way using black-box elements corresponding to $n$-cycles or
$(n-1)$-cycles and $3$-cycles. 

Thus the 
elements we wish to construct correspond to  $m$-cycles with $m \in
\{2, 3, n-1, 
n\}$. Such elements $g$ satisfy 
the equation $g^m=1$. For
$t=2$ or $t=3$ the proportion of  $t$-cycles in $S_n$ or $A_n$ is very
small.  Typically the  algorithms  construct  elements  $g$ whose
cycle  structure 
consists of  a $t$-cycle and a  single additional non-trivial 
cycle  of length $r$, where $r$ is
not  divisible  by $t$ (as in Table~1), and $g$ satisfies 
$g^{tr} =1$  and $g^r \not= 1$.
For such an element we  construct a $t$-cycle
by forming the  power $g^r$. 
We note  that, in the case where  $n\equiv 1\pmod{6}$ and $n
\ge 13$, we may utilise elements in  both cases 9 and 10 of Table~1 to
construct a  3-cycle (although if  $n=13$, some additional  care is  needed).

It  turns out that  most ele\-ments  in $S_n$  or $A_n$  satisfying an
equation of the form $g^m=1$ are $m$-cycles if $m\in \{n, n-1\}$. Also
most  elements for which  $g^{tr}=1$ and  $g^r \not=  1$ consist  of a
$t$-cycle and an $r$-cycle if $t=2$ or $3$.  The crucial probabilistic
result  underpinning the algorithms  in \cite{Bealsetal03}  shows that
the conditional probability that a random element $g$ has one of these
desired cycle types, given that $g$ satisfies an appropriate equation,
is $1-o(1)$  for large $n$, and at  least $1/180$ for all  $n$.  In an
algorithmic  context  this  means  that  a  random  element  of  $S_n$
satisfying  one of these  equations has  a good  chance of  having the
desired  cycle structure  and the  lower bounds  in \cite{Bealsetal03}
were sufficient  for the purpose  of estimating the complexity  of the
algorithms.   However, for  an efficient  practical  implementation of
these algorithms a more realistic  lower bound is desirable, since the
lower bound  is reciprocally  proportional to an  upper bound  for the
number of random  elements that need to be  tested.  In particular, if
the  algorithm is called  with an  input group  $G$ not  isomorphic to
$A_n$ or $S_n$  then the number of random  elements considered will be
equal to the upper bound.

\section{On Divisors of Integers}\label{sect:divisors}

In this section we cite some results from Number Theory which we
require throughout. In particular we investigate properties
of divisors of a given integer and sums of powers of these divisors.
For a positive integer $n$ let $d(n)$ denote the number of divisors of
$n$. Niven et al. \cite[pp. 395-396]{NivenZuckermanetal91} prove the
following result: 

\begin{lemma}\label{lem:nrdiv}
For every $\delta > 0$ there is a constant $c_{\delta}$ such that 
$d(n) \le c_\delta n^{\delta}$ 
for all $n \in \N.$
In particular, we may take
$c_{1/2} = \sqrt{3}$ and $c_{1/3} =(1536/35)^{1/3}\sim 3.53.$
\end{lemma}

In the following  two lemmas we revisit and refine  the proof given in
 \cite{NivenZuckermanetal91}
to obtain certain constants  $c_0 < c_{1/3}$ such that $d(n) \le c_0
 n^{1/3}$ for most integers $n$. 

\begin{lemma}\label{lem:oddnrdiv}
Let $n\in \N$.

\begin{tabular}{ll}
{\rm(a)}& If $n$  is odd  then $d(n) \le 4 (3/35)^{1/3} n^{1/3}< 1.764\,
n^{1/3}.$\\  
{\rm(b)}& If $n$ is not divisible by $9$ then
$d(n) \le (16/105^{1/3}) n^{1/3}$\\
&$ < 3.392\, n^{1/3}$.\\
{\rm(c)}&
If $n$ is odd and  not 
divisible by $9$, then 
$d(n) \le (8/105^{1/3}) n^{1/3}$\\
&$ <1.696\, n^{1/3}.$\\ 
\end{tabular}
\end{lemma} 

\begin{proof}
Let $n = \prod_p p^{\alpha_p}$, where the product is over all odd
primes $p$ and each $\alpha_p \ge 0$. Then,
arguing as in   
\cite[pp. 395-396]{NivenZuckermanetal91}, we have that   
$$\frac{d(n)}{n^{1/3}} = \prod_{p\mid n} f_p(\alpha_p),$$ where
  $f_p(\alpha) = \frac{\alpha+1}{p^{\alpha/3}}$.  
It is shown in \cite[pp. 395-396]{NivenZuckermanetal91} 
that the function $f_p(\alpha)$ of an integer 
variable $\alpha$ attains its maximum at  
 $\alpha_0(p) = \lfloor
\frac{1}{p^{1/3}-1} \rfloor.$   In particular, if $p >
8$ then
$\alpha_0(p) =  0$   and  thus   $f_p(\alpha)  \le  1, $ for all $\alpha$.
Also $\alpha_0(5) = \alpha_0(7) = 1$,\ $\alpha_0(3) = 2$,
and $\alpha_0(2) = 3.$ 

If $n$ is odd, then $\alpha_0(2)=0$ and hence by
taking $c = f_3(2) f_5(1) f_7(1) = 4 (3/35)^{(1/3)}$ 
we have that
$d(n) \le c n^{1/3}.$ This proves (a).

If $n$ is not divisible by $9$, then $\alpha_0(3) \le 1$ and hence by
taking $c = f_2(3) f_3(1) f_5(1) f_7(1) = 16/(105)^{(1/3)}$ 
we have that $d(n) \le c n^{1/3}.$ This proves (b).

Finally, if $n$ is odd and not divisible by $9$ then by
taking $c = f_3(1) f_5(1) f_7(1) = \frac{8}{105^{1/3}}$ 
we have that $d(n) \le c n^{1/3}$. This proves (c).
\end{proof}

\bigskip

\begin{lemma}\label{lem:bignrdiv}
Let $c_0 = (\frac{768}{35})^{1/3}\sim 2.8$. Then  
for $n \in \N,$ either $d(n) \le c_0\, n^{1/3}$, or
 $n = 2^{a} 3^{b} 5^{c} 7^{d} m$, where
$1 \le a \le 6$, $0 \le b \le 4$, 
$0 \le c \le 2 $, $ 0\le d \le 1$,  and $m \in \{ 1, 11, 13\}$. In
particular, if $n > 11,793,600$ then $d(n) \le c_0 n^{1/3}.$
\end{lemma}

\begin{proof}
Let $n = \prod_p p^{\alpha_p}$, and let $f_p, \alpha_p$ and
$\alpha_0(p)$ be as in the proof of Lemma~\ref{lem:oddnrdiv}.
Note that 
 $\alpha_0(2) = 3,$  $\alpha_0(3) = 2,$ and 
$\alpha_0(5) = \alpha_0(7) = 1.$
Let $c(n) = \prod_{p\mid n}
f_p(\alpha_p),$  so that $d(n) = c(n) n^{1/3}.$
Let $n_0 = 2^7 \cdot 3^2\cdot 5 \cdot 7.$
Then $$c(n_0) =
f_2(7) f_3(2) f_5(1) f_7(1) =
\frac{8}{2^{7/3}}\frac{3}{3^{2/3}}\frac{2}{5^{1/3}}\frac{2}{7^{1/3}} 
= c_0,$$ and if $c(n) \le c_0$ then we obtain $d(n) \le c_0 n^{1/3}.$ 

Now suppose $n$ is such that $c_0 < c(n).$ 
Write $n$ as
 $n = 2^{a} 3^{b} 5^{c} 7^{d} m$, where
$\gcd(m, 2\cdot  3\cdot 5\cdot 7 ) = 1$. Now 
 $c(n) =  f_2(a) f_3(b) f_5(c) f_7(d)c(m),$ 
 and $c(m) \le 1$ with $c(m) < 1$ if $m > 1.$ The condition $c_0 <
 c(n)$, and our knowledge of the maximum values for the $f_p(\alpha)$
 give 
\begin{eqnarray*}
c_0 = f_2(7) f_3(2) f_5(1) f_7(1) &<& f_2(a) f_3(b) f_5(c) f_7(d)c(m)\\
& \le &  
f_2(a) f_3(2) f_5(1) f_7(1)
\end{eqnarray*}
and hence 
$f_2(7) < f_2(a)$ which implies  $a \le 6.$

Next, it is easy to see that $f_p(\alpha) < f_q(\alpha)$ for $p > q$
and any positive integer $\alpha$. Also for
any $p > 8$ the function  $f_p(\alpha)$ is decreasing for $\alpha \ge
0.$

Suppose that one of the following holds:
\begin{enumerate}
\item[(i)] $m$ is divisible by some prime $p$ with $p \ge 17$;
\item[(ii)] $m$ is divisible by $11^2$ or $13^2$ or $11\cdot 13.$
\end{enumerate}

We show that in either case $f_2(3)c(m) < f_2(7).$ 
In case (i),\\ $f_2(3) c(m) \le f_2(3) f_{17}(1) < f_2(7)$. In case (ii),
since $f_{13}(2) < f_{11}(2)$, we have 
either $f_2(3) c(m) \le f_2(3) f_{11}(2) < f_2(7)$ or 
\[f_2(3) c(m) \le f_2(3) f_{11}(1) f_{13}(1) < f_2(7).\] 
Thus if (i) or (ii) holds, then $f_2(3) c(m) < f_2(7)$ and we have
\begin{eqnarray*}
c_0 &= & f_2(7) f_3(2) f_5(1) f_7(1) \\
&<& c(n) 
=  f_2(a)f_3(b) f_5(c) f_7(d) c(m) \\
&\le&f_2(3) f_3(2) f_5(1) f_7(1) c(m) \\
&<& f_2(7) f_3(2) f_5(1) f_7(1) = c_0,
\end{eqnarray*}
which is a contradiction. Thus $m \in \{1, 11, 13\}.$
Also $f_2(3) f_7(2) < f_2(7) f_7(1)$, and a similar argument yields
$d \le 1;$ and $f_2(3)f_5(3) < f_2(7)f_5(1)$ and we obtain $c\le 2;$
and $f_2(3) f_3(5) < f_2(7)f_3(2)$ and so $b \le 4.$ Finally
$f_3(4) f_5(2)f_7(1) = 2 \left(\frac{5}{21}\right)^{1/3} < c_0$ and it
follows that $a \ge 1.$ 
\end{proof}

The following lemma yields some  elementary approximations.

\begin{lemma}\label{lem:sumdiv}Let $n, a, b $ be positive integers with $a
\le b \le n$.  Let $D$ denote the set 
of all divisors $d$ of $n$ for which $a \le d \le b.$ Then
 $$\sum_{d\in D} (d-1)(d-2) \le (b-1)(b-2)  + nb -na.$$
\end{lemma}

\begin{proof}
A divisor $d \in D$ is of the form $d = n/t$ for some integer $t$.  
As $a\le d \le b$ it follows that $ n/b \le  t \le n/a.$ Thus 
\begin{eqnarray*}
\sum_{d\in D} (d-1) (d-2) & \le& 
\sum_{t=n/b}^{n/a} (\frac{n}{t}-1)(\frac{n}{t}-2) \\
& \le&   (b-1)(b-2) + \int^{n/a}_{n/b}
\frac{n^2}{t^2} {\rm d} t\\
& = &  (b-1)(b-2)  + nb -na.
\end{eqnarray*}
\end{proof}

\section{Estimating Proportions of elements}\label{sec:proportions}

Let $m$ and $n$ be positive integers with $m \ge n-1$.
We estimate the proportion $P(n,m)$ of elements in the 
symmetric group  $S_n$   whose order divides $m$.  
Note that the order $|g|$ of a permutation $g \in S_n$
divides $m$ if and only if the length of each cycle of $g$
divides $m$.  
 Thus $P(n,m)$ is the proportion of elements in $S_n$ all of
whose cycle lengths divide $m$. 
As indicated in the introduction, we obtain estimates for
proportions of elements in $S_n$ whose order divides $m$ 
in various ways. We begin by defining different proportions of
elements which play a key role in our analysis.

\begin{notation}\label{notation}
Let
$P^{(1)}(n,m)$ denote the proportion of elements $g\in S_n$ of order
dividing $m$ for which $1,2,3$ lie in the same $g$-cycle, 
let
$P^{(2)}(n,m)$ denote the proportion of elements $g\in S_n$ of order
dividing $m$ for which $1,2,3$ lie in exactly two $g$-cycles and 
let
$P^{(3)}(n,m)$ denote the proportion of elements $g\in S_n$ of order
dividing $m$ for which $1,2,3$ lie in three different $g$-cycles.

\end{notation}

Note that 
\begin{equation}\label{eq:p}
P(n,m) = P^{(1)}(n,m) + P^{(2)}(n,m) + P^{(3)}(n,m)
\end{equation}
and 
 that by convention
we take $P(r,m) = 1$ if $r \le 0.$

We begin by deriving expressions for $P^{(i)}(n,m)$ for $i=1,2,3$.

\begin{lemma}\label{lem:theps}
Let $n$ and $m$ be positive integers with $m \ge n-1$.
Then the following all hold, where we take $P(0,m) = 1.$
\begin{enumerate}
{\openup 7pt
\item[(a)]
${\displaystyle P^{(1)}(n,m) = \frac{(n-3)!}{n!} 
\sum_{{\stackrel{\scriptstyle d \mid m}{3 \le d\le
      n}}}{(d-1)(d-2)}P(n-d,m).}$\\
\item[(b)]
${\displaystyle P^{(2)}(n,m) =
  \frac{3(n-3)!}{n!}\sum_{\stackrel{\scriptstyle d_1, d_2     \mid m
    }{2\le  d_2,\, d_1+d_2\le n}}  (d_2-1)P(n-d_1-d_2,m). }$
\item[(c)]
${\displaystyle P^{(3)}(n,m) = \frac{(n-3)!}{n!}
  \sum_{\stackrel{\scriptstyle d_1,d_2,d_3\mid  m}{d_1+d_2+d_3  \le n}}  
P(n-d_1-d_2 -d_3,m). }$
\item[(d)]
${\displaystyle P(n,m) = \frac{1}{n} \sum_{d |m,\, 1\leq d\leq n} P(n-d,m)}.$

}
\end{enumerate}
\end{lemma}

\begin{proof}
We first compute $P^{(1)}(n,m)$, the proportion of those permutations in $S_n$
for which  the points
 $1, 2, 3$ are contained in one $g$-cycle, $C$ say, of 
length $d$ with $d \mid m$ and $3\le d.$ Also $d \le n$ since $g \in
S_n.$

We can choose the remainder of the  support set of $C$ in $\binom{n-3}
{d-3}$ ways and then 
the cycle $C$ in $(d-1)!$ ways. 
The rest of the
permutation $g$ can be chosen in $P(n-d,m)(n-d)!$ ways. Thus,
for a given 
$d$, the number of such elements is
$(n-3)!(d-1)(d-2)P(n-d,m)$. We obtain
the proportion $P^{(1)}(n,m)$ by summing over all  divisors 
$d$ of $m$ which are at most $n$, and dividing the sum by $n!$, that is
$$P^{(1)}(n,m) = \frac{(n-3)!}{n!}  
\sum_{\stackrel{\scriptstyle d \mid m}{3 \le d\le n}}{(d-1)(d-2)P(n-d,m)}.$$ 
Hence part
(a) follows.  Parts (b) and (c) are derived in a similar fashion. For
a detailed proof of a very similar result see the proof of Lemma~2.2
of 
\cite{NiemeyerPraeger05b}.

Also part (d) follows by enumerating the 
elements $g\in S_n$ of order dividing $m$ according to the length $d$ of the 
$g$-cycle containing the point~1. 
\end{proof}

We now  prove Theorem~\ref{first}.

\begin{proof}
Let $\gamma$ denote $\gamma(m)$ as in Definition~\ref{def:gamma}.
The result is immediate if
$\frac{1}{n} + \frac{\gamma m}{n^2} \ge 1$ since $P(n,m) \le 1$, so we
may assume that $1 > \frac{1}{n} + \frac{\gamma m}{n^2} > \frac{\gamma
  m}{n^2},$ whence $n > \sqrt{\gamma m}.$  In particular $n \ge 4$
and if $n = 4$ then $m = 3$. However $
{\displaystyle P(4,3) = \frac{9}{24} = \frac{3}{8} < \frac{1}{4} +
  \frac{\gamma(3) 3}{16}}$. Thus  we may assume that  $n \ge 5.$
Let $D$ denote the set of all divisors of $m$
which are at most $n$.

Using the fact that $P(t,m)\leq 1$ for $t < n$ 
in Lemma~\ref{lem:theps}(1) we obtain 
\begin{equation}\label{eq:a}
P^{(1)}(n,m) \le\frac{(n-3)!}{n!} 
\sum_{\stackrel{\scriptstyle d \in D}{d\ge 3}}{(d-1)(d-2)}.
\end{equation}
By applying Lemma~\ref{lem:sumdiv}  we obtain, if $m \ge n$: 
\begin{eqnarray*} 
P^{(1)}(n,m) & \le & \frac{1}{n(n-1)(n-2)} \left( (n-1)(n-2) + mn - 3m \right)\\
& = &  \frac{1}{n} + \frac{m(n-3)}{n(n-1)(n-2)}\\
& < &  \frac{1}{n} + \frac{m}{n^2}\\
\end{eqnarray*}
and similarly
\begin{eqnarray*} 
P^{(1)}(n,n-1) & \le & \frac{ (n-2)(n-3) + n(n-4)}{n(n-1)(n-3)}\\
& = &  \frac{1}{n} + \frac{n(n-4)-(n-3)}{n(n-1)(n-3)}\\
& < &  \frac{1}{n} + \frac{n-1}{n^2}\\
& = &  \frac{1}{n} + \frac{m}{n^2}.\\
\end{eqnarray*}
Now let $D_2 = \{ (d_1, d_2) \colon d_1, d_2 \in D,\, 2 \le d_2,\, d_1 +
d_2 \le n \}.$ Then,
using the fact that $P(t,m)\leq 1$ for $t < n$ 
in Lemma~\ref{lem:theps}(2) we obtain 
\begin{eqnarray}
P^{(2)}(n,m) & \le  & \frac{3(n-3)!}{n!}\sum_{(d_1, d_2) \in D_2} (d_2-1)
\label{eq:b}\\ 
& =  & \frac{3}{n(n-1)}\sum_{(d_1, d_2) \in D_2} \frac{d_2-1}{n-2}.\nonumber
\end{eqnarray}
Since $d_1 + d_2 \le n$ and  $1 \le d_1$ it follows that
$d_2-1 \le n - d_1 - 1 \le n-2.$ Set $c(m) = d(m)/m^{1/3},$ where
$d(m)$ is the number of divisors of $m$. Then
\begin{eqnarray*}
P^{(2)}(n,m) 
& \le  & \frac{3}{n(n-1)}\sum_{d_1, d_2 \in D} 1 \\
& = & \frac{3m^{2/3}}{n(n-1)} c(m)^2.
\end{eqnarray*}

Using the fact that $P(t,m)\leq 1$ for $t < n$ 
in Lemma~\ref{lem:theps}(3) we obtain 
\begin{eqnarray}
P^{(3)}(n,m) &\le& \frac{(n-3)!}{n!} \sum_{d_1,d_2,d_3\in D} 1\label{eq:c}\\
&=&  \frac{c(m)^3m}{n(n-1)(n-2)}.\nonumber 
\end{eqnarray}

Now using the inequality 
$n > \sqrt{\gamma m}$ 
in the upper bounds for the $P_i(n,m)$ gives
\begin{eqnarray*}
P(n,m)  &\le  & 
\frac{1}{n} + \frac{m}{n^2} +  \frac{3c(m)^2 m^{2/3}}{n(n-1)} 
+ \frac{c(m)^3m}{n(n-1)(n-2)} \\
& < &  \frac{1}{n} \!+ \!\frac{m}{n^2} \left( 1\! +\!
\frac{3c(m)^2\sqrt{\gamma m}}{m^{1/3}(\sqrt{\gamma m}-1)} +
\frac{c(m)^3\sqrt{\gamma m}}{(\sqrt{\gamma m}-1)(\sqrt{\gamma m}-2)}
 \right).
\end{eqnarray*}
Consider the function $$f(m,c) =
\frac{3c^2}{m^{1/3}}\frac{\sqrt{\gamma m} }{(\sqrt{\gamma m}-1)} + 
\frac{c^3\sqrt{\gamma m}}{(\sqrt{\gamma m}-1)(\sqrt{\gamma m}-2)}.$$ 
If  $c(m) \le c$ then 
$$P(n,m) < \frac{1}{n} + \frac{m}{n^2}\left( f(m,c)+1\right).$$

Thus, to prove Theorem~\ref{first} for any given value of $m$ (and for 
all $n$  with $n-1\leq m$), it is sufficient to prove that $f(m,c)\leq 
\gamma(m) -1$ for some $c\geq c(m)$. This is the way we shall  obtain  our 
result for large $m$.

Recall that  $c(m) =  d(m)/m^{1/3}$ and as in
Lemma~\ref{lem:bignrdiv}, we set $c_0 = (768/35)^{1/3}.$ 
It is easy to see that, for fixed $c$, the function 
$f(m,c)$ is strictly decreasing as $m$ increases over any interval on
which $\gamma = \gamma(m)$ is constant. Moreover, we
can check that $f(19020, c_0) \le 1.$ Thus if both 
 $m \ge 19020$ and    $c(m) \le c_0,$ then
 $f(m,c_0) \le 1 = \gamma - 1$
and hence $f(m,c) \le \gamma - 1$, so the theorem is
proved in this case.

The remaining values are
 all $m < 19020$ and 
those  $m \ge 19020$ for which $c(m) > c_0$.
Note that Lemma~\ref{lem:bignrdiv} identified explicitly 
a finite set of integers that contains all integers
$m$ such that
$c(m) > c_0$. For each of these remaining $m$ we need to consider all
 $n$ such that $\sqrt{\gamma m} < n \le m+1$. 
We define
$${S}(n,m) = \sum_{\stackrel{\scriptstyle d\mid m}{3\le d \le
    n}}{(d-1)(d-2)} +  3\sum_{\stackrel{\scriptstyle d_1, d_2 \mid
    m}{2 \le d_2,\, d_1 + d_2     \le n}}  
(d_2-1)  + \sum_{ \stackrel{\scriptstyle d_1,d_2,d_3\mid
    m}{d_1+d_2+d_3 \le n}}  1.$$ Then 
 Equations~($\ref{eq:a}$),($\ref{eq:b}$) and ($\ref{eq:c}$) imply 
that $P(n,m) \le \frac{(n-3)!}{n!} S(n,m).$ Thus to prove $P(n,m) \le
\frac{1}{n} + \frac{\gamma m}{n^2}$ it is sufficient to prove that
$S(n,m) \le (n-1)(n-2) (1 + \frac{\gamma m}{n})$
(for each given $m$ and all $n\leq m+1$).

Next we define  
$$\hat{S}(n,m) = \sum_{\stackrel{\scriptstyle d\mid m}{3\le d \le
    n}}{(d-1)(d-2)} +  
     3\sum_{(d_1, d_2) \in D_2} (d_2-1) 
+  \sum_{ (d_1,d_2,d_3 ) \in D_3} 1,$$ 
where 
$D_2 = \{ (d_1, d_2 ) \colon  d_i \le n,\, d_i\! \mid\! m,\, 2 \le
    d_2, \, d_1 + 
    d_2 \le m \}$ and  
$D_3 = \{ (d_1, d_2, d_3 ) \colon  d_i \le n,\, d_i \mid m,\, d_1 + d_2 +
    d_3 \le m\}.$  
Then 
$S(n,m) \le \hat{S}(n,m)$, so it suffices to prove that
\begin{equation}\label{eq:star}
\hat{S}(n,m) \le (n-1)(n-2) ( 1 + \frac{\gamma m}{n})
\end{equation}
(for each given $m$ and all $n\leq m+1$).
Note that for fixed $m$ the right hand side of Inequality
($\ref{eq:star}$) is increasing in $n$. Thus suppose $d$ divides $m$ and 
$\hat{S}(d,m) \le (d-1)(d-2) ( 1 + \frac{\gamma m}{d}).$
If $d < n \le m+1$ and there is no divisor of $m$ in the interval
$(d,n]$ then
every divisor $d_i$ of $m$ that satisfies $d_i \le n$ also satisfies
$d_i \le d,$ and hence 
\begin{eqnarray*}
\hat{S}(n,m) = \hat{S}(d,m) &\le& (d-1)(d-2) ( 1 + \frac{\gamma m}{d}
)\\
& \le & (n-1)(n-2) ( 1 + \frac{\gamma m}{d}).
\end{eqnarray*}
Hence, if for all $d$ with $d\mid m$ and $d > \sqrt{\gamma m}$ we have
$\hat{S}(d,m) \le (d-1)(d-2) ( 1 + \frac{\gamma m}{d})$ then
$S(n,m) \le (n-1)(n-2) (1 + \frac{\gamma m}{n})$ for all $n \le m + 1.$

For   all  $m   \le  19020$   and  for all   $m  \ge   19020$   for  which
$c(m)>c_0=(768/35)^{1/3}$  is possible as given by
Lemma~\ref{lem:bignrdiv}, we 
tested in {\sf GAP} \cite{GAP4} whe\-ther,  for  all divisors $d$ of
$m$, the inequality 
$ \hat{S}(d,m) \le (d-1)(d-2) (  1 + \frac{\gamma m}{d})$ holds, where
$\gamma= \gamma(m)$ is as in Definition~$\ref{def:gamma}$.  This was
the case for  all values of 
$m$ we tested, except for $m = 72$ and $m =120.$  In these two
exceptional cases 
we proved directly that $S(n,m) \le (n-1)(n-2)(1 + \gamma m/n)$ holds
for all $n$ 
with $\sqrt{\gamma m} \le n \le m+1.$ 
 Thus the theorem is proved.
\end{proof}

\section{Proof of Theorem~\ref{melbounds_sn}}\label{sec:longcycle}
 
In this section we give a full proof of the probability bounds in
Theorem~\ref{melbounds_sn}. 

We repeatedly use the following arithmetic fact which holds
for all $x,y >0$ 
\begin{equation}\label{eq:fact}
\frac{x}{x+y} > 1 - \frac{y}{x}.
\end{equation}
Our next result refines the upper bound for $P(n,m)$ in Theorem~\ref{first}
for the special case where $m$ is $n$ or $n-1$. It deals with cases
$1,4,5$ of Table~1.

\begin{theorem}\label{the:prac}
 Suppose that $n,m$
are positive integers  and $m \in \{n-1, n\}$ such that one of the
cases $1, 4,$ or $5$ of Table~$1$ holds with $m = r.$
Let $g$ be  a uniformly distributed  random ele\-ment from $S_n$ 
(for case $1$),  or $A_n$ (for cases $4$ and $5$). 
Let $A, B$ denote the events that $g$ is an  $m$-cycle, or $g$
has order dividing $m$, respectively.
Let $\gamma(m)$ be defined as in Definition~$\ref{def:gamma}$.
 Then
\begin{enumerate}
\item[(a)] ${\displaystyle 
P(n,m) \le  \frac{1}{m} + \frac{d(m)}{n^2}   (2 + 4\gamma(m)) 
},$ and
\item[(b)] 
{\openup 3pt
${\displaystyle
P(A\mid B) \ge 1 - \frac{(2+4\gamma(m))d(m)}{n}}$. 
Moreover,  for all cases, \\
$P(A\mid B) \ge 1/2$ except for the cases in case~$1$ of Table~$1$ where}
$m = n$ divides $24$, and in these exceptional cases
$P(A\mid B) \ge 2/7.$
\end{enumerate}
\end{theorem}

\noindent
\begin{proof}
(a) \quad
By Lemma 3.2~(d), $P(n,m)=\frac{1}{n}\sum_{d|m, d\leq n}P(n-d,m)$. The 
values of $d$ in the summation satisfy either $d=m$ or $d\leq m/2$. Since 
$P(0,m)=P(1,m)=1$, it follows that

$$P(n,m) = \frac{1}{m } + 
\frac{1}{n} \sum_{\stackrel{\scriptstyle d\mid m}{d \le m/2}} P(n-d,m).$$

As $d \le m/2 \le n/2$ we have that $n-d \ge   n/2.$ Thus by 
 Theorem~\ref{first} we obtain that
\begin{eqnarray*} 
P(n,m) &\le &\frac{1}{m} +\frac{1}{n}  
\sum_{d \mid m,\, d\le m/2}\left( \frac{1}{n-d} +
\frac{\gamma(m) m}{(n-d)^2}\right) \\ 
 &\le &\frac{1}{m} +\frac{1}{n}  
\sum_{d \mid m,\, d\le m/2}\left( \frac{2}{n} +
\frac{4\gamma(m) m}{n^2}\right) \\ 
 &\le & \frac{1}{m} + \frac{1}{n}  \sum_{d \mid m,\, d\le m/2} \frac{2
  + 4\gamma(m)}{n}\\  
 &\le & \frac{1}{m} + \frac{d(m)}{n^2}   (2 + 4\gamma(m)).\\ 
\end{eqnarray*}
This proves part (a).

\bigskip
\noindent
(b)\quad  
As $A \subseteq B$ the conditional probability $P(A\mid B)$
satisfies $P(A\mid B) = P(A\cap B)/P(B) = P(A)/P(B)$.
For case 1 of Table~1, where $g$ is chosen from $S_n$, we have $P(B)=P(n,m)$.
For cases $4$ and $5$ of Table~1 where $g$ is chosen from $A_n$, we have
$P(B) \le 2P(n,m)$ since the number of elements $g \in A_n$
satisfying  $g^m = 1$ is at most equal to the number of such elements
in $S_n$.

Further, in the case of $S_n$ (case~1 of Table~1), 
 $P(A)$ is the proportion of $m$-cycles in $S_n$, which is
$1/m \ge 1/n$, since $m = n$ or $m=n-1.$ 
In the case of $A_n$ (cases~4 and 5 of Table~1),
$m$ is odd,  and $P(A) = 2/m \le 2/n.$ Hence in all three cases, using part 
(a), $P(A\mid B)$  satisfies
\begin{equation*}
 P(A\mid B) \ge 
\frac{\frac{1}{m}} {P(n,m)} \ge
\frac{1} {1 + \frac{d(m) ( 2 + 4 \gamma(m))}{n} }. 
\end{equation*}
By (\ref{eq:fact}) it follows that 
$$ P(A\mid B) \ge 
1 - \frac{d(m)}{n} ( 2 + 4 \gamma(m) ).$$

For $n  \ge 362$  we have $m \ge n-1 > 360$ and hence $\gamma(m)=2$ and by
Lemma~\ref{lem:nrdiv} we get
$$ P(A\mid B) \ge 
\frac{1} {1 + \frac{10 (1536/35)^{1/3}}{n^{2/3}}} \ge 
\frac{1} {1 + \frac{10 (1536/35)^{1/3}}{362^{2/3}}} > 1/2.$$

We used the approximation $P(A\mid B) \ge 
(\frac{1}{m})/P(n,m)$ and computed $P(n,m)$ precisely (using
Lemma~\ref{lem:theps}(d)) ,   
to verify by computation  in {\sf  GAP} 
that the conditional probability $P(A\mid B) \ge 1/2$ for all values
of $n, m$ as in the statement. For the remaining cases, that is those
in case~1 of Table~1 where $m=n$ divides $24$, we computed that the
lower
bound $(nP(n,m))^{-1}$ for $P(A\mid B)$ is greater than 2/7.
\end{proof}

Now we turn to  determining the conditional probability that 
 an element $g\in S_n$ has cycle structure $2^1r^1$, given that
$g$ has order dividing $2r$ and $|g^r| =2$, where $r \in \{n-2, n-3 \}$
and $r$ is odd.

\begin{remark}\label{lem:divmttw} 
Let $n$ and $r$  be positive integers such that $n \ge 7$, 
 $r \in \{n-2, n-3\}$, and $r$ is odd.
If $d$ is a divisor of $2r$ with $d \le n$
then either $d = r$, or $d = 2r/3$, or $d \le 2r/5.$
\end{remark} 

\begin{theorem}\label{the:melpractw}
Let $n$, $r$ and {\sc CycType}  be as in
case~$2$ or $3$ of Table~$1$, with $n\geq8$.
Let $g$ be  a uniformly distributed  random ele\-ment from $S_n$, 
and let $A, B$ denote the events that $g$ has cycle structure $2^1r^1$, or 
$g$ has order dividing $2r$ and $|g^r|=2$, respectively.
Let $\gamma(2r)$  be as in
Definition~$\ref{def:gamma}$.
Then 
\begin{enumerate}
\item[$(a)$] 
$P(B) \le {\displaystyle
\frac{1}{2r} +
\frac{1}{n^2}(3 + 18\gamma(2r) ) +
 \frac{d(2r)}{n^2}(\frac{5}{3} + \frac{50 \gamma(2r)}{9})}$, and
\item[$(b)$] $\displaystyle{
 P(A \mid B) > 1 - 
\frac{1}{n}(6 + 36\gamma(2r) ) -
 \frac{d(2r)}{n}(\frac{10}{3} + \frac{100\gamma(2r)}{9}).
}$\\
For $n \not\in \{11,17,18\}$, $P(A\mid B)$ is  at least $1/3$ while
 for $n \in \{11,17,18\}$ it is   at least $1/4.$ 
\end{enumerate}

\end{theorem}

\begin{proof}
\quad (a)\quad
If $g\in S_n$ has order dividing $2r$, then the length $d$ of any
$g$-cycle divides $2r$ and $d\le n.$ By Remark~\ref{lem:divmttw}, 
 $d = r$ or $d = 2r/3$ or $d \le 2r/5.$ 
We divide $B$ into two disjoint events $B_1$ and $B_2$, where $B_1$
is the event that $g$ contains an $r$-cycle and $B_2$ is the event
that it does not. Then $P(B) = P(B_1) + P(B_2)$. 

If $g$ has a cycle of length $r$ and
if $|g^{r}|=2$ then $|g|$ divides $2r$ and $g$ has cycle type
$2^1r^1$. Hence $B_1$ is equal to the event $A$, and the probability
that a random element of $S_n$ has cycle structure $2^1r^1$ is
$$P(B_1) = \frac{\binom{n}{r} (r-1)! \, (n-r)!}{2(n!)}  = \frac{1}{2r}.$$  

Let  $D'$ denote the set of all divisors of $2r$ which are at most
$2r/3$. 
Let $P'(n,2r)$ denote the proportion of elements of $S_n$ all of whose
cycle lengths lie in $D'$.
Then by Remark~\ref{lem:divmttw}, for any element $g \in S_n$ such that $B_2$ holds, the $g$-cycle
containing the point 1 has length $d$ for some $d \in D'$. For a given
$d \in D'$, we estimate the number of possible $g$ as follows.

We have $\binom{n-1}{d-1}(d-1)!$ 
choices of $d$-cycles containing 1 and at most
  $P'(n-d,2r) (n-d)!$
choices for the rest of the permutation. Summing over all
divisors $d \in D'$ yields 
$$P(B_2) \le \frac{1}{n}\sum_{d \in D'}P'(n-d,2r).$$

By Theorem~\ref{first}
 we obtain 
$$\frac{1}{n}\sum_{d \in D'}P'(n-d,2r) 
\le
 \frac{1}{n}  \sum_{d \in D'} 
\left(
\frac{1}{n-d} + \frac{2r \gamma(2r) }{(n-d)^2} 
\right).$$

 If $d = 2r/3$ then
$n-d = n -\frac{2r}{3} > \frac{n}{3}$, while if $d \le 2r/5$
then 
$n-d \ge n -\frac{2r}{5}>\frac{3n}{5}.$ Hence
we obtain 
\begin{eqnarray*}
P(B_2) &< & 
\frac{1}{n}(\frac{3}{n} + \frac{18\gamma(2r)}{n}) +
 \frac{1}{n}  \sum_{d \in D', d\not= 2r/3} 
\left(
\frac{5}{3n} + \frac{50\gamma(2r)}{9n} 
\right) \\
&\le &
\frac{1}{n^2}(3 + 18\gamma(2r) ) +
 \frac{d(2r)}{n^2}(\frac{5}{3} + \frac{50 \gamma(2r)}{9}).
\end{eqnarray*}
Adding this bound to $P(B_1)$ yields part (a).

\medskip
(b) \quad 
Since $A \subseteq B$ it follows that $P(A\mid B) = P(A)/P(B).$
We showed in the proof of part (a) that $P(A) = 1/(2r).$
Thus, using part (a) we obtain
\begin{eqnarray*}
P(A\mid B) &\ge & 
\frac{\frac{1}{2r}}{\frac{1}{2r} + 
\frac{1}{n^2}(3 + 18\gamma(2r) ) +
 \frac{d(2r)}{n^2}(\frac{5}{3} + \frac{50 \gamma(2r)}{9})}\\
&> & 
\frac{1}{ 1 +  
\frac{1}{n}(6 + 36\gamma(2r) ) +
 \frac{d(2r)}{n}(\frac{10}{3} + \frac{100 \gamma(2r)}{9})}.\label{eq:x}
\end{eqnarray*}
Finally, by  Inequality~(\ref{eq:fact}),
$$P(A\mid B) > 
1 - \frac{1}{n}(6 + 36\gamma(2r) ) -
 \frac{d(2r)}{n}(\frac{10}{3} + \frac{100 \gamma(2r)}{9}).$$
This proves the first assertion of part (b).

Since $r$  is odd we obtain by Lemma~\ref{lem:oddnrdiv}
that $d(2r) = 2 d(r) \le 8 \left(\frac{3}{35}\right)^{1/3} n^{1/3}$.
Thus (\ref{eq:x}) yields
$$
P(A\mid B) >
\frac{1}{ 1 +  
\frac{1}{n}(6 + 36\gamma(2r) ) +
 \frac{8\left(\frac{3}{35}\right)^{1/3}}{n^{2/3}}
(\frac{10}{3} + \frac{100 \gamma(2r)}{9})}.
$$
If $n\ge 360$, then  $\gamma(2r) = 2$ and
this lower bound on the
conditional probability is at least $0.3335 > 1/3.$ For smaller values 
of $n$ we proceed as follows. Note that 
 $P(A\mid B) = P(A)/P(B)$ and  $P(B)$
is the proportion of elements of order dividing $2r$ but not of order
dividing $r$. Thus $P(A\mid B) = 
\frac{\frac{1}{2r}} {P(n,2r)-P(n,r)}$.
By computation  in {\sf  GAP} 
we verified, by computing $P(n,2r)$ and
$P(n,r)$  precisely  (using Lemma~\ref{lem:theps}(d)),  
that the conditional probability $P(A\mid B)$ satisfies the lower
bounds given in the statement.
\end{proof}

\noindent
\emph{Proof of Theorem}~\ref{melbounds_sn} (a) and (b). \quad
(a) It follows immediately from Theorem~\ref{the:prac}(b), and
Lemma~\ref{lem:nrdiv}, that for the events $A$ and $B$ as defined in 
Theorem~\ref{the:prac},
$P(A\mid B) \ge 1 -
(\frac{1536}{35})^{(1/3)}\frac{(2+4\gamma(m))}{n^{2/3}}
> \frac{8 + 15\gamma(n)}{n^{2/3}}$ (since $\gamma(m) \le \gamma(n)$).
The absolute lower bounds for this probability were proved in
Theorem~\ref{the:prac}(b).

(b)  Since $r$ is odd, $d(2r) = 2 d(r)$, and so by
Lemma~\ref{lem:oddnrdiv}, $\gamma(2r) \le 8\left(
\frac{3}{35}\right)^{1/3}n^{1/3}$. Then, by 
Theorem~\ref{the:melpractw}, with $A$ and $B$ the events defined
there, 
$P(A\mid B) \ge 1 - \frac{
(6 + 36\gamma(2r) ) +  8 \left(\frac{3}{35}\right)^{1/3} 
(\frac{10}{3} + \frac{100 \gamma(2r)}{9})}{n^{2/3}} > 1 - \frac{18 +
    76\gamma(2r)}{n^{2/3}}.$
The absolute lower bounds for this probability were proved in
Theorem ~\ref{the:melpractw}(b). \hfill
$\fbox{\hphantom{.}}$

\bigskip

It remains to 
determine the conditional probability that
an element $g$ of $S_n$ or $A_n$ has 
cycle structure {\sc CycType} as in one of the cases 6-10 of Table~1,
given that $g$ has order dividing $3r$ with $n, r$ and
$|g^r|$  as in that case of  Table~1. We will deduce
Theorem~\ref{melbounds_sn}(c)  from the following result.

\begin{theorem}\label{the:melpracthree}
Let $n, r$  and {\sc CycType}  be  as in one of the cases $6-10$ of
Table~$1,$ and let $n \ge 8$.
Let $g$ be a uniformly distributed  random element $g$ in $S_n$ (for 
cases $6-9$) or $A_n$ (for case $10$). Let $A$ denote
the event that $g$ has cycle type {\sc CycType} as given in the 
relevant case of Table~$1$, and let $B$ denote the event that $g^{3r}=1$
and $|g^r|=3$. 
Let $\gamma = \gamma(3r)$ $($see Definition~$\ref{def:gamma})$
so that $\gamma$ satisfies
$${\displaystyle
\gamma = \begin{cases} 2 & {if\ } n \ge 124, \\
2.5 & {if\ } 26 \le n \le 123, \\
3.345 & {if\ } n \le 25. 
	   \end{cases}}$$
Then
\begin{enumerate}
\item[$(a)$]  
${\displaystyle P(B) \le 
\begin{cases} 
 \frac{c}{3r} +  \frac{7+39\gamma}{n^2} +
 \frac{d(3r)}{n^2} \left( \frac{20+75\gamma}{16}\right), & \mbox{for\
 cases\ } 6-9\\
 \frac{1}{3r} +  \frac{8+96\gamma}{n^2} +
 \frac{d(3r)}{n^2} \left( \frac{10+75\gamma}{2}\right), & \mbox{for\
 case\ } 10\\
\end{cases}}
$
\noindent
$$\mbox{where\ }\quad \begin{array}{ll}
c = \begin{cases} 1/2 & {for\ case\ } 8,\\
                    1/3 & {for\ case\ } 9,\\
                   1 &   {for\ cases\ } 6,7.
      \end{cases} 
\end{array}$$
\item[$(b)$] 
For cases $6-9$,
$$P(A \mid B) \ge  a\left(1 -
 \frac{3(7+39\gamma)}{cn} -
\frac{d(3r)}{n}\frac{3(20 + 75\gamma)}{16c}\right),
$$ 
where $a = 1$ for cases $6-8$ and $a=1/2$ for case $9,$ and $c$ is as
in $(a)$. Also, if in  
the events $A$ and $B$ the random element is restricted to lie in
$A_n$, then $P(A\mid B)$ is unchanged from its value in $S_n$.

\noindent
For case $10$,
$$P(A \mid B) \ge  1 -
 \frac{3(8+96\gamma)}{n} -
\frac{d(3r)}{n}\frac{3(10 + 75\gamma)}{2}.
$$ 
Moreover, the lower bounds on $P := P(A \mid B)$ given in Table~$2$ hold.
\end{enumerate} 
\end{theorem}

The following technical result, used in the proof of
Theorem~\ref{the:melpracthree} can be proved using similar techniques 
to those used in the proof of Theorem~\ref{the:melpractw}.

\begin{lemma}\label{lem:divmt} 
Let $n, r$ be as in one of the cases $6-10$ of Table~$1$ with $n \ge
8$. Let $d$ be a divisor
of
 $3r$ such that $d \le n$. Then for cases $6-9$ one of 
$d= r$, or $d \le r/5$, or $d = 3r/y$ with $y \in 
\{ 5, 7, 11, 13 \},$  and for case $10$, one of 
$d= r$, or $d = 3r/4$, or $d \le 3r/5$,
or $(n,r,d)=(13,8,12)$.
\end{lemma} 

\begin{proof} 
Suppose $d$  is a divisor of  $3r$ and $d \le n,$ say $d =
3r/c.$ Then $3(n-6) \le 3r \le cn$.
 We claim 
that $c\geq 2.$ Suppose to the contrary that $c=1.$
 Then $n\leq 9$, so $n$ is 8 
or 9 and in either case $r=5$. Thus $9\geq d=3r=15$, which is a 
contradiction. So $c\geq 2$.

In cases 6-9, $r \equiv \pm 1 \pmod{6}$ and hence $c$ is
odd and $c \ge 3;$ the values $c = 3,5,7,11,13$ give possibilities
listed. Since $r \equiv \pm 1 \pmod{6}$ we know $3r \not\equiv
0\pmod{9}$ and hence $c \not= 9.$ This leaves $c \ge 15$ which implies
$d \le r/5.$

Finally in  case 10,  $r=n-5 \equiv  2 \pmod{6}$ and  $c \geq  2.$ The
values $c=3, 4$ and $c \geq 5$ give  the possibilities listed. 
The remaining value
$c = 2$  corresponds to $d= 3(n-5)\le n$ which implies  $n = 13$,
$r=8,$ and $d=12.$
\end{proof}

\emph{Proof of Theorem}~\ref{the:melpracthree}. \quad
(a)\quad
If $g\in S_n$ has order dividing $3r$, then the length $d$ of any
$g$-cycle divides $3r$ and $d\le n.$ 
By Lemma~\ref{lem:divmt}, 
either $d\leq r,$ or we have case 10 with 
$(n,r,d)=(13,8,12)$. However in this exceptional case, the element $g$ 
would be a $12$-cycle, contradicting the fact that $g$ is an even 
permutation. Thus in all cases $d\leq r$.
Note that $r > 1$ as $n \ge 8.$
We divide $B$ into two disjoint events $B_1$ and $B_2$, where $B_1$
is the event that $g$ contains an $r$-cycle and $B_2$ is the event
that it does not. Then $P(B) = P(B_1) + P(B_2)$. 

If $g$ has a cycle of length $r$ and $|g^{r}| = 3$ then 
$|g|$ divides $3r$ and $g$ has cycle type $3^1 r^1$ (for cases 6-8),
 $3^1 r^1$  or $3^2 r^1$  (for case 9), or $2^13^1r^1$ (for case 10).
Hence $P(B_1)$ is the proportion  of such elements in $S_n$ (for cases
6-9) or $A_n$ (for case 10), namely
\begin{equation}
P(B_1) =  \begin{cases}
P(A) = \frac{1}{3r} & \text{for\ cases\ } 6, 7 \mbox{\ or\ } 10\\
P(A) = \frac{1}{6r} & \text{for\ case\ } 8\\
2P(A) = \frac{1}{9r} & \text{for\ case\ } 9.\\
  \end{cases}
\label{cases:pb1}
\end{equation}

Let  $D'$ denote the set of all divisors of $3r$ which are less than
$r$. Then by Lemma~\ref{lem:divmt}, for cases 6-9,  $D' \subseteq D_1
\cup D_2$, where
$D_1 = \{ d \in D' \mid d \le 
r/5\}$ and $D_2 = D' \cap \{ 3r/y \mid y \in \{5, 7, 11, 13\}\}$, and
for case 10, $D'$ is the set of all divisors $d$ of $3r$ with $d \le
3r/4.$ Then $P(B_2)$ is the proportion of elements $g$ in $S_n$ (for
cases 6-9) or $A_n$ (for case 10), all of whose cycle lengths lie in
$D'$ and for which $|g^r| = 3.$ Let $P'(n,3r)$ denote the proportion
of elements in $S_n$ of order dividing $3r$ and not containing an
$r$-cycle. Then $P(B_2) \le P'(n,3r)$ (for cases 6-9) and $P(B_2) \le
2P'(n,3r)$ (for case 10). 

We shall now estimate $P'(n,3r)$. This is the proportion of elements
of $S_n$ with all cycle lengths in $D'$. Considering the elements $g
\in S_n$ with all cycle lengths in $D'$ according to the length $d$ of
the $g$-cycle containing the point 1, we have
$$P'(n,3r) = \frac{1}{n}\sum_{d \in D'}P'(n-d,3r).$$ 
Suppose first that $n, r$ are as in one of the cases 6-9, so $D'
\subseteq D_1 \cup D_2$. If $d \in D_1$, then $n-d \ge n-r/5 > 4n/5.$
Note also that $n-d \le n-1 \le r + 5 \le 3r$ (since $r > 1$ and $r$
is odd). Then by 
Theorem~\ref{first}, 
\begin{eqnarray}
\frac{1}{n}\sum_{d \in D_1}P'(n-d,3r) 
&\le&
 \frac{1}{n}  \sum_{d \in D_1} \left(
\frac{1}{n-d} + \frac{\gamma 3r}{(n-d)^2} 
\right) \label{eq:first}\\
&<&
 \frac{d(3r)}{n}  \left(
\frac{5}{4n} + \frac{75\gamma}{16n} 
\right)\nonumber \\
&=&
 \frac{d(3r)}{n^2}  (\frac{5}{4} + \frac{75}{16} \gamma),\nonumber
\end{eqnarray}
where $\gamma$ is as in the statement.
Next we consider the divisors in $D_2.$ In this case  $d = 3r/y$
 where $y \in \{5, 7, 11, 13\}$. Then
 $n-d = n-3r/y > (y-3)n/y.$ Hence applying
 Theorem~\ref{first} we have, with $\gamma$ as in the statement,
\begin{eqnarray*}
\lefteqn{\frac{1}{n}\sum_{d \in D_2}P'(n-d,3r)  }\\ & \le & \frac{1}{n}
\sum_{d \in D_2} \left( \frac{1}{n-d} + \frac{\gamma 3r}{(n-d)^2} \right) \\
&< & \frac{1}{n^2}\sum_{y\in \{5, 7, 11, 13\}}\left( \frac{y}{(y-3)}
+ \frac{\gamma 3y^2}{(y-3)^2} \right)\label{eq:ys}\\ 
& =  &\frac{277}{40n^2} +  \frac{61887 \gamma }{1600n^2}\\
& <  &  \frac{7+39\gamma}{n^2}.\\
\end{eqnarray*}

Thus, for cases 6-9, 
\begin{equation}\label{eq:melprac1}
P(B_2) \le P' (n,3r) \le
 \frac{7+39\gamma}{n^2} +
 \frac{d(3r)}{n^2}  \frac{20+75\gamma}{16}
.\end{equation}
This proves 
(a) for cases 6-9, since $P(B) \le P(B_1) + P'(n,3r).$

Now consider case 10. Here $r = n-5$ and for $d \in D'$ either $d =
3r/4$ and $n-d > n/4$, or $d \le 3r/5$ and $n-d > 2n/5.$ Also
$n-d \le n-1 \le r + 4 \le 3r$ and hence by Theorem~\ref{first}, 
\begin{eqnarray*}
\frac{1}{n}\sum_{d \in D'}P'(n-d,3r) 
&\le&
 \frac{1}{n}  \sum_{d \in D'} \left(
\frac{1}{n-d} + \frac{\gamma 3r}{(n-d)^2} 
\right) \\
&<&
\frac{1}{n^2} (4 + 48 \gamma) +
 \frac{d(3r)}{n^2}  (\frac{5}{2} + \frac{75\gamma}{4}). \\
\end{eqnarray*}

Hence
\begin{eqnarray}
P(B_2) \le 2P'(n,3r) &< &
 2\frac{4 + 48 \gamma}{n^2} + 2d(3r)
 \frac{(10 + 75\gamma)}{4n^2}\nonumber \\
&=&
 \frac{8 + 96 \gamma}{n^2} + d(3r)
 \frac{(10 + 75\gamma)}{2n^2}.\label{eq:lastline} 
\end{eqnarray}
This proves $(a)$ for case 10.

\medskip
(b) \quad 
For cases 6-9 let $P'(A), P'(B)$ denote the probabilities of the
events $A$ and $B$ respectively in the case where the
random element is chosen only from $A_n$. Since $3r$ is odd, all
elements of order dividing $3r$ lie in $A_n$ and so $P'(A) = 2P(A),$
and $P'(B) = 2P(B)$. Consequently, the conditional probability
$P(A\mid B)$ remains the same for these cases if we restrict the
elements to lie in $A_n$.  Hence we only need to consider random
elements
in $S_n$ for cases $6-9$.

First we consider the case where $A$ is the event that a
uniformly distributed random element $g$ in $S_n$ (for cases 6-9) or
$A_n$ (for case 10) has cycle type $3^1r^1$ (for cases 6-9) or
$2^13^1r^1$ (for case 10). In this case
we see from the discussion above that 
$P(A) = a P(B_1)$, where $a = 1$ for cases $6, 7, 8$ and $10$, and
$a=1/2$ for case 9. Note $P(B_1)$ is given by $(\ref{cases:pb1})$,
so $P(B_1) = c/(3r) > c/(3n)$ with $c$ as in part (b) for cases
$6-9$ and $c=1$ for case 10. Since $A \subseteq B$, the conditional
probability satisfies 
\begin{eqnarray}
P(A \mid B) &=&  \frac{P(A)}{P(B)}  =  
\frac{aP(B_1)}{P(B_1)+P(B_2)} \label{eq:prac0}\\
\frac{1}{1+\frac{3n}{c}P(B_2)}&>& a\left(1 -
\frac{P(B_2)}{P(B_1)}\right). \nonumber 
\end{eqnarray}
For cases 6-9
the general assertions of (c) follow immediately from this, the 
fact that $r \le n,$ and from (\ref{eq:melprac1}).
For case 10, the assertion follows from this, the fact that $r\le n$ and from 
(\ref{eq:lastline}).

Finally we have to show that the conditional probabilities are at
least $1/3$ for all cases except $(n,r) = (31,25)$.

For cases 6-9, $r\equiv \pm 1\pmod{6}$ so that $3r$ is odd and not
divisible by $9$. Thus by Lemma~\ref{lem:oddnrdiv}, 
$$d(3r) < \frac{8}{(105)^{1/3}} (3r)^{1/3} < \frac{8}{(35)^{1/3}}
n^{1/3},$$ so by (\ref{eq:melprac1}), 
$$P(B_2) \le \frac{7 + 39\gamma}{n^2} + \frac{8}{(35)^{1/3}
  n^{5/3}}\left( \frac{20 + 75\gamma}{16}\right).$$ 
For $n \ge 124$ we have $\gamma = 2$ and so 
\begin{equation}\label{eq:prac1}
P(B_2) \le 
\frac{85}{n^2} + \frac{85}{(35)^{1/3}
  n^{5/3}}.
\end{equation}

By Equations ($\ref{eq:prac0})$ and ($\ref{eq:prac1}),$ 
the  conditional probability $P(A\mid B)$ is at least 
$$
P(A\mid B) \ge 
\frac{aP(B_1)}{
P(B_1) + 
 \frac{85}{n^2} +
 \frac{85}{(35)^{1/3}n^{5/3}} }.
$$
In cases 6 and 7 we obtain
$$
P(A\mid B) \ge 
\frac{\frac{1}{3r}}{
\frac{1}{3r} +
 \frac{85}{n^2} +
 \frac{85}{(35)^{1/3}n^{5/3}} }
>
\frac{1}{
1 +
 \frac{255}{n} +
 \frac{255}{(35)^{1/3}n^{2/3}} }.
$$
For $n \ge 420$ this is at least $1/3.$

\noindent
In case 8 we obtain
$$
P(A\mid B) \ge 
\frac{\frac{1}{6r}}{
\frac{1}{6r} +
 \frac{85}{n^2} +
 \frac{85}{(35)^{1/3}n^{5/3}}  }
>
\frac{1}{
1 +
 \frac{510}{n} +
 \frac{510}{(35)^{1/3}n^{2/3}} }.
$$
For $n \ge 1050$ this is at least $1/3.$

\noindent
In case 9 we obtain
$$
P(A\mid B) \ge 
\frac{\frac{1}{18r}}{
\frac{1}{9r} +
 \frac{85}{n^2} +
 \frac{85}{(35)^{1/3}n^{5/3}}  }
>
\frac{1}{
2 +
 \frac{1530}{n} +
 \frac{1530}{(35)^{1/3}n^{2/3}} }.
$$
For $n \ge 12400$ this is at least $1/3.$

In case 10 
since $r = n-5 \equiv 2\pmod{6}$ it follows 
that $3r$ is not divisible by $9$. By Lemma~\ref{lem:oddnrdiv},
$d(3r) \le \frac{16}{(105)^{1/3}} (3r)^{1/3}  <
\frac{16}{(35)^{1/3}} n^{1/3}$ and hence by (\ref{eq:lastline})
$$P(B_2) \le 
\frac{8 + 96 \gamma}{n^2} + \frac{8}{(35)^{1/3} n^{5/3}} 
(10 + 75\gamma).$$
If $n \ge 124 $ we have $\gamma =2 $ and hence
\begin{equation*}
P(B_2)  <  \frac{200}{n^2} + \frac{1280}{(35)^{1/3} n^{5/3}}.
\end{equation*}
Therefore, by ($\ref{eq:prac0}$)
$$
P(A\mid B) \ge \frac{\frac{1}{3r}}{
\frac{1}{3r} +
 \frac{200}{n^2} +
 \frac{1280}{(35)^{1/3}n^{5/3}}  }
>
\frac{1}{
1 + 
 \frac{600}{n} + 
 \frac{3840}{(35)^{1/3}n^{2/3}}  }.
$$
For $n\ge 14700$ this is at least $1/3.$

For the remaining values of $n$ we showed that the statement holds
by computation in {\sf GAP}. For sufficiently large $n$ we computed
all divisors of $3r$ and evaluated (\ref{eq:first}) directly for cases
$6-9$, to obtain a better upper bound for $P(B_2)$ than that in
($\ref{eq:melprac1}$), and similarly for case 10. For some
values of $n$ this was not sufficient to show that $P(A\mid B) \ge
1/3.$ For these values of $n$ we computed the proportions precisely to
obtain the lower bounds given in the statement. \hfill
$\fbox{\hphantom{.}}$

\bigskip
Finally we prove Theorem~\ref{melbounds_sn}(c)
using the results proved in Theorem~\ref{the:melpracthree}.
\smallskip

\noindent
\emph{Proof of Theorem}~\ref{melbounds_sn} (c). \quad
Note first that the 
absolute lower bounds for the conditional probability in Table~$2$
are proved in Theorem~\ref{the:melpracthree}(b).

For cases $6-8$, $3r$ is 
odd and not divisible  by 9, so by Lemma~\ref{lem:oddnrdiv}, $d(3r) \leq 
\frac{8}{(105)^{1/3}} (3r)^{1/3} 
< \frac{8}{(35)^{1/3}} r^{1/3} 
$. Also, for these cases, in 
Theorem~\ref{the:melpracthree}(b)
we have $a=1$ and $c\geq 1/2$, so
\begin{eqnarray*}
P(A \mid B) 
&\ge&  \left(1 - \frac{2\cdot3(7+39\gamma)}{n} -
\frac{1}{(35)^{(1/3)}n^{2/3}}\frac{2\cdot 3(20 + 75\gamma)}{2}\right)\\
&\ge&  \left(1 - \frac{42+234\gamma + (60 + 150\gamma)/(35)^{1/3}}
{n^{2/3}}\right)\\
&\ge&  1 - \frac{61+303\gamma}{n^{2/3}},
\end{eqnarray*}
which is greater than $1-\frac{98+839\gamma}{n^{2/3}}.$

For case 9 again we have $d(3r) \le 8/(35)^{(1/3)} n^{1/3}$ but this
time the parameters $a$ and $c$ of Theorem~\ref{the:melpracthree}(b)
have the values  $a = 1/2$ and $c=1/3.$ Thus
\begin{eqnarray*}
P(A \mid B) 
&\ge&  \frac{1}{2}\left(1 - \frac{3\cdot3(7+39\gamma)}{n} -
\frac{1}{(35)^{(1/3)}n^{2/3}}\frac{3\cdot 3(20 + 75\gamma)}{2}\right)\\
&=& \frac{1}{2} - \left(\frac{63 +351\gamma}{2n} +
\frac{1}{(35)^{(1/3)}n^{2/3}}\frac{180 + 675\gamma}{4}\right)\\
&\ge& \frac{1}{2} - \frac{46 +228\gamma}{n^{2/3}}.
\end{eqnarray*}

Finally, for case 10, $r \equiv 2\pmod{6}$ so
$3r$ is not divisible by $9$.
Hence by Lemma~\ref{lem:oddnrdiv}(b),
 $d(3r) \le 16/105^{(1/3)} (3r)^{1/3} < 16/35^{(1/3)} n^{1/3}$.
Thus
\begin{eqnarray*}
P(A \mid B) 
&\ge&  \left(1 - \frac{3(8+96\gamma)}{n} -
\frac{8}{(35)^{(1/3)}n^{2/3}}3(10 + 75\gamma)\right)\\
&\ge&  1 - \frac{98+839\gamma}{n^{2/3}}.
\end{eqnarray*}
\hfill$\fbox{\hphantom{.}}$

\subsection*{Acknowledgements}
The authors acknowledge the support of ARC Discovery Project DP0557587.

{\footnotesize
\providecommand{\bysame}{\leavevmode\hbox to3em{\hrulefill}\thinspace}

}

\end{document}